# BOLZANO'S CONJECTURE: MEASURING THE NUMEROSITY OF INFINITE SETS.

## JULIAN JACK

The mere fact, therefore, that two sets A and B are so related that every member a of A corresponds by a fixed rule to some member b of B in such wise that the set of these couples (a+b) contains every member of A or B once and only once, never justifies us, we now see, in inferring *the equality of the two sets, in the event of them being infinite,* with respect to the multiplicity of their members ……….On the contrary, and in spite of their entering symmetrically into the above relation with one another, the two sets can still stand in a relation of inequality, in the sense that the one is found to be a whole and the other a part of that whole. [6, page 98].

**ABSTRACT**
Bolzano and Cantor were the first mathematicians to make significant attempts to measure the size (numerosity) of different infinite collections. They differed in their methodological approaches, with Cantor's prevailing. This led to the foundation of the theory of sets as well as Cantor's transfinite arithmetic. This paper argues that Bolzano's conjecture is correct and that Euclid's principle, 'that the whole is greater than a part', should be considered as a necessary condition for the quantification of infinite sets (rather than bijection). Cantor had concluded that the rational and the algebraic numbers were of the same size as the natural numbers, whilst, in contrast, the real numbers were a larger set. Using Cantor's methods it is shown in this paper that the rational numbers are of larger size than the natural numbers, thus showing that bijection is not a reliable measure of the size of infinite sets. It is also concluded, using mathematical induction, that different 'countably' infinite sets can have various different sizes. The implication for theorems using bijection as a measure of size is then briefly discussed. There already exist new methods of measuring numerosity, based on Euclid's principle, which may develop a consistent system of infinite arithmetic.

**INTRODUCTION**
In the nineteenth century two mathematicians sought to develop an arithmetic of infinite quantities. Both Bolzano and Cantor were in favour of an actual or completed infinity and they also agreed that it should be possible to compare and count the sizes of different infinite collections. They differed in their methodological approach to counting. Cantor used a specific form of pairing (one-to-one correspondence) called bijection to decide whether two collections were of equal size (or cardinality). One problem was that it was already known that this assumption, for infinite collections, was incompatible with Euclid's principle, 'that the whole is greater than a part' [16]. Bolzano, by contrast, thought Euclid's principle should be given precedence over bijection. Bolzano's last, posthumous, publication on the subject was in 1851. Cantor, who was born in 1845, subsequently produced a sequence of publications [e.g. 7, 8] that led to his version of counting to be widely accepted and this was reinforced by his larger body of work which formed the basis for axiomatic set theory. Amongst the striking results arising from Cantor's assumptions was the conclusion that there were different sizes of infinite collections, with the sizes being quantified by the (infinite) cardinal numbers. The smallest infinite cardinal number, $\aleph_0$, denoted the size of the collections with a countable number of elements (amongst others, including the natural numbers, the rational numbers and the algebraic numbers), with larger cardinal numbers ($\aleph_1, \aleph_2, \aleph_3 \ldots$) then being generated by a power law (Cantor's power theorem). Despite the acceptance of his work there has been persistent uncertainty and controversy [19], the principal issue being the assignment of a cardinal number to the real numbers (the continuum hypothesis), since they had been shown, by Cantor's other method of size comparison, which was qualitative rather than quantitative, to be a larger collection than the natural numbers.

Cantor gave two proofs [7, 8], that the real numbers were more numerous that the natural numbers, both of which were indirect proofs (argument by contradiction) of which the simpler, later, proof [8] relied on the diagonal method.

In contrast to Cantor's form of infinite arithmetic, Euclid's principle would result in a different hierarchy of infinite sizes with, for example, the natural numbers, the rational numbers, the algebraic numbers and the real numbers all being of successively greater size. This paper shows, using the diagonal method in the same way as Cantor, that subsets of the rational numbers, of a particular order, can be greater in size than the natural numbers. Since the method of bijection, in Cantor's system, concludes that the natural numbers and the rational numbers are of the same size, the system allows the generation of contradictions i.e. it is inconsistent. Explanations are offered for the reason that the rational numbers, unlike the real numbers, have to have a particular order and also for why a bijective function is not a sufficient condition to measure infinite size.

**PROOF OF CONTRADICTION**

Consider the set made up of decimal numbers, where each number is defined by $\sum_{i=1}^{n} 10^{-i}$, where $n$ goes from 1 to $\infty$; i.e. an ordered list of decimals made up exclusively of the digits 0 or 1, starting with the decimal equivalent of $\frac{1}{10}$ and terminating with the equivalent of $\frac{1}{9}$. This set will be called the original set. In order to give a one-to-one correspondence between the members of this set and the natural numbers, the last decimal on the list (0.1111111...... or $0.\dot{1}$) has to be paired with the first natural number and the first decimal (0.1 or 0.1000000......) with the second natural number, the second decimal (0.11 or 0.1100000......) with the third natural number, etc. Note that this method of pairing is that suggested for showing that sets of order type $\omega$ and $\omega+1$ have the same cardinality (e.g. Figure 1, p.617 in J. Bagaria, IV.22 in [9]). In Cantor's terminology, this set has an ordinality of a member of the second ordinal number class ($\omega$ +1) and a cardinality of the natural numbers.

Diagonalization on this set (in which the $n^{th}$ digit of the $n^{th}$ number is changed) will always produce a new number that is not a member of the set, even if the diagonal rule is to replace 0 with 1 and *vice versa*. However, using this diagonal rule, it is possible to generate a diagonal number that is a member of an extended version of this set. It is simply a matter of choosing to add as elements to the start of the set, one or more (m) distinct rational numbers which each have a 0 in the $n^{th}$ decimal place. An example, for m=5, is shown, where the decimal versions of 1/100, 1/1000, 1/10000, 1/100000 and 1/1000000 are put in the first five (vertical) ordinal positions. This extended set retains the cardinality of the natural numbers on Cantor's bijective assumption.

```
        1 2 3 4 5 6 7 8 9 ......n.......      (horizontal ordinal number, of digits)
1    0.0 1 0 0 0 0 0 0 0..............
2    0.0 0 1 0 0 0 0 0 0..............
3    0.0 0 0 1 0 0 0 0 0..............
4    0.0 0 0 0 1 0 0 0 0..............
5    0.0 0 0 0 0 1 0 0 0..............
6    0.1 0 0 0 0 0 0 0 0..............
7    0.1 1 0 0 0 0 0 0 0..............
8    0.1 1 1 0 0 0 0 0 0..............
9    0.1 1 1 1 0 0 0 0 0 ..............
.        ..............................
```

| n | ............................ |
| . | ............................ |

The diagonal number will be $0.\dot{1}$, which is in ordinal position $\omega$ on the list. The key assumption of the diagonal method, when applied to numbers between 0 and 1 expressed in their decimal form, is that since the number of digits in each decimal number and the number of numbers in the set are both bijectable with the natural numbers, a new number formed by changing a single digit of each number (as long as the change occurs for a digit at a different ordinal position in each number) cannot be a member of the set, since all the numbers in the set must have been changed, if the number of numbers and the number of digits in each number are the same size. If the diagonal number is a member of the set then the number of numbers must be larger than the number of digits in each number. Thus, on Cantor's assumptions, the above extended set of rational numbers is larger in size than the natural numbers.

A second example of this deduction comes from considering a set of decimal numbers in which the numbers in the odd vertical ordinal position are defined by $\frac{1}{99 \times 10^{(n-1)}}$ and those in the even ordinal position by $\sum_{i=1}^{\frac{n}{2}} 10^{-i}$ (where n goes from 1 to ∞). This leads to a set with the following appearance:

```
      1 2 3 4 5 6 7 8 9 . . . . . n . . . . . . .    (horizontal ordinal number)
1     0.0 1 0 1 0 1 0 1 0 . . . . . .  . . . . . . .
2     0.1 0 0 0 0 0 0 0 0 . . . . . .  . . . . . . .
3     0.0 0 0 1 0 1 0 1 0 . . . . . .  . . . . . . .
4     0.1 1 0 0 0 0 0 0 0 . . . . . .  . . . . . . .
5     0.0 0 0 0 0 1 0 1 0 . . . . . .  . . . . . . .
6     0.1 1 1 0 0 0 0 0 0 . . . . . .  . . . . . . .
7     0.0 0 0 0 0 0 0 1 0 . . . . . .  . . . . . . .
.     . . . . . . . . . . . . . . . . . . . . . . .
.     . . . . . . . . . . . . . . . . . . . . . . .
.     . . . . . . . . . . . . . . . . . . . . . . .
n     . . . . . . . . . . . . . . . . . . . . . . .
.     . . . . . . . . . . . . . . . . . . . . . . .
.     . . . . . . . . . . . . . . . . . . . . . . .
.     . . . . . . . . . . . . . . . . . . . . . . .
```

followed by $0.\dot{0}$ in the $\omega$ position and $0.\dot{1}$ in the $\omega+1$ position (i.e. this set is of order type ($\omega$ +2)). Once again, a diagonalization using the same rule yields a number, this time in the $\omega+1$ position.

The general rule for examples of this proof of different sizes of any countably infinite set is that, for the members of the set, a list is generated with the order of the second ordinal number class and a member in ordinal position $\omega$ or beyond is selected as the 'target', with the members in the first part of the ordinal listing then being selected, along with the diagonal rule, so that the 'target' is generated by the diagonalization. This applies, not only to decimal numbers, but also to sets made up of other members – such as sets containing the elements m and w, as in Cantor's original 1891 proof [8] – where the diagonal rule simply changes m to w and *vice versa*.

**CAN THE DIFFERENCE IN SIZE BE 'QUANTIFIED'?**
Two specific examples of sets whose elements are a particular order of rational decimal numbers have been given above, in which it is concluded that the diagonal number is a member of the set. For each of them, there is a constructive proof, using mathematical induction, to deduce how the size of the set compares with the size of the natural numbers. Let L(n) specify the n$^{th}$ number in the list and d(n) the n$^{th}$ digit of each number. The simplest way to construct a diagonal number is to start with d(1) of L(1), then d(2) of L(2) etc. Let D(n) specify the diagonal number as it is being formed iteratively. Applying mathematical induction to each set in turn:

**First example.**
Base case: D(1)=L(6).
Induction step: If D(n)=L(n+5), both equal $\sum_{i=1}^{n} 10^{-i}$.
> To form D(n+1), note that d(n+1) in L(n+1) is 0. Therefore, on the diagonal rule, d(n+1) in D(n+1) is 1. In other words $10^{-(n+1)}$ must be added to D(n) to form D(n+1).
> D(n+1)= $\sum_{i=1}^{n} 10^{-i} + 10^{-(n+1)} = \sum_{i=1}^{n+1} 10^{-i}$ =L(n+6).

Conclusion: For all n ≥ 1, D(n) is a member of the list, five positions further down the list. On completion (i.e. diagonalization over all the natural numbers), D(∞) is the same number as the last member of the list and thus the set has a size of 5 more than the natural numbers (and of the original set). This conclusion will apply when the m of 5 is replaced by a different finite number of insertions at the start of the set, with the size differing by m.

**Second example.**
Base case: D(1)= L(2).
Induction step: If D(n)= L(2n), then both equal $\sum_{i=1}^{n} 10^{-i}$.
> To form D(n+1), note that d(n+1) in L(n+1), for both odd and even values of n, is 0. Therefore d(n+1) in D(n+1) is 1 and $10^{-(n+1)}$ has to be added to D(n) to form D(n+1).
> D(n+1)= $\sum_{i=1}^{n} 10^{-i} + 10^{-(n+1)} = \sum_{i=1}^{n+1} 10^{-i}$ = L(2n+2).

Conclusion: For all n (≥ 1), D(n) occurs at a position of 2n on the list. Therefore, the list is twice the size of the natural numbers (and of the original set).

**DISCUSSION**
The key conclusion of this paper is that selected subsets of the rational numbers (Q) can be shown, by the diagonal method, to be of a larger size than the natural numbers i.e. in Cantor's terminology, uncountable. An assertion to the contrary is made by Abbott [1, Exercise **1.6.3**]. "Rebut the following about the proof of **Theorem 1.6.1**: every rational number has a decimal expansion, so we could apply this same argument to show that the set of rational numbers between 0 and 1 is uncountable. However, because we know that any subset of Q must be countable the proof of **Theorem 1.6.1** must be flawed." **Theorem 1.6.1.** of Abbott states that the reals in the open interval 0 to 1 are uncountable (by the diagonal method). The solution offered by Abbott is as follows: "if we imitate the proof to try to show that Q is uncountable, we can construct a real number x in the same way. This x will again fail to be in the range of our function f, *but there is no reason to expect x to be rational.* The decimal expansions for rational numbers either terminate or repeat, and this will not be true of the constructed x". Whilst this assertion is usually true, with random order of the numbers, this paper shows that there is an infinite number of counterexamples, when one selects subsets of rational numbers and puts them in a particular order. To understand this latter condition, it may be helpful to consider the diagonal method, when applied to a finite list of numbers each with a finite number of (significant) digits. Let us restrict our consideration to numbers with only two possible digits (so that there is only one appropriate diagonal rule), with the number of digits and the minimal number of numbers being n. There are only two sizes of the list of numbers which are certain, by diagonalization, to yield an unambiguous result: when the list is of size n, the result is a diagonal product not in the list and when the list is of size $2^n$ (the largest it can be without repeating a number in the list) the diagonal product must be a member of the list. The reason for such a large range of uncertainty for sizes between (n+1) and ($2^n - 1$) is that it will depend on both the order and selection of the numbers (out of the possible) in the list.

As an illustration consider a finite set of decimal numbers between 0 and 1 with 5 significant digits, restricting the possible digits at each place to 0 or 1. There are 32 possible numbers. This set will be named rea. The numbers can be listed in 32! ($\cong 2.63 \times 10^{35}$) different orders. Diagonalization (with the rule of 0 to 1 and vice versa) will cover the first five numbers in the order and the number generated will be somewhere in the remaining 27 numbers further down the order, whichever order is selected, because all possible numbers are in the set.

Define a subset of rea: rat is those numbers in rea which terminate either in three 0's or three 1's. There are 8 of these numbers and they can be ordered in 40,320 different ways. If one selects the 'target' number as the diagonal product 0.11111 then the first five numbers in the order must have 0 in the digital position equivalent to their place in the order and the 'target' number must be located in one of the last three positions in the order. These restrictions lead to the calculation that only 720 of the possible 40,320 orders will lead to a diagonal product which is a member of the set (rat). This illustrates how, for finite arithmetic, diagonalization is not, in general, a definitive probe of the size of the set, unless that set includes all possible elements. Of course, by selecting any of the other seven possible 'target' numbers, the number of possible orders which yield a definitive result from diagonalization increases, but it is still a small proportion of the total number of possible orders.

This analysis of diagonalization in a finite system provides insight into the reason why, in the infinite case, order (including order type) is important unless the set on which the diagonalization is performed contains all the possible numbers in the chosen interval (i.e. the relevant real numbers). The nature of the ordering and selection that is required in the case of infinite sets of rational numbers for the diagonal product to be a member of the list is given in the final section of the Proof of Contradiction.

This paper has shown that bijection is not a reliable measure of the size of infinite sets, although it is for finite sets. It is clear, in retrospect, why this arises. The principle of one-to-one is a method of pairing that requires a mapping function so that one element of one set is exclusively paired with only one element of the second set. This requirement always leads to the same conclusion about size when both sets have a finite number of elements, despite the fact that there are a large number of different possible mappings. The same result is not true for the comparison of two infinite sets (for an historical account, see Mancosu [14]). For example, if a finite collection with m elements is added to an infinite collection there exist mappings that can 'prove' equality or leave unpaired from 1 to m of the elements of the extended collection. Mancosu [14] gives the example of two infinite sets, the even numbers and the natural numbers (starting with 1). The pairing with the mapping function 2n↔n gives equality whereas the mapping function 2n↔2n leaves all the odd numbers unpaired (suggesting the natural numbers are twice the size of the even numbers). Cantor decided to ignore this ambiguity and elected to make the existence of at least one mapping function, which concluded equality (bijection), as the necessary and sufficient condition for judging the size of infinite sets. This meant he had to have another method to prove the inequality of two infinite sets and both of his proofs [7, 8] were indirect (argument by contradiction) which, in classical logic, requires a consistent system to be reliable. To make bijection reliable for infinite sets, the sufficient condition would be to establish that all possible one-to-ones reached the same conclusion about size, as for finite sets. The conclusion that bijection does not guarantee equal size has considerable implications, particularly for proofs using the diagonal method as well as any indirect proofs, which may not be reliable. Thus Cantor's power theorem does not exclude the possibility of sizes intermediate between $\aleph_0$ and $\aleph_1$, etc. His proof about the existence and construction of transcendental numbers (see Gray [10]) is not reliable because, depending on the order of the algebraic numbers, the diagonal product may itself be an algebraic number. The implications for the application of the diagonalization lemma, as used in the first incompleteness theorem of Gödel, also needs to be carefully explored.

It may be noted that in order to reach the conclusions about the 'quantification' of the size of the sets, using mathematical induction, it was assumed that quantification over all the natural numbers was possible. In simple terms this means that the possibility of a completed or actual infinity was assumed. Detailed discussions of this from different points of view may be found in articles by Linnebo and Shapiro [12] and by Zenkin [21]. A more important point is that this method of 'quantification' is not a general method, since it is dependent on the order of the elements in the sets. There is therefore a need to develop a way of measuring the numerosity of infinite sets which is not sensitive to the order of the elements.

The results presented here strongly reinforce Mancosu's argument [13, 14, 15] that there is no inevitability to Cantor's assumption that, for infinite sets, the existence of a one-to-one mapping function showing equinumerosity should be a sufficient condition for measuring size. Mancosu drew attention to the work of others (Katz [11] and Benci and his colleagues e.g. [ 3, 4, 5]) who have been developing new methods of measuring the size of infinite sets. These methods both incorporate Euclid's principle as a necessary condition for the measurement. One important result, from Benci and DiMasso [3], is that their measurement of numerosity is not only consistent but independent of the axioms of set theory as well, so that a consistent system for measuring numerosity could be conjoined with an axiomatic theory (such as Peano's arithmetic) and the result would be consistent, as long as the axiomatic theory is consistent.

Recently, Bellomo and Massas [2] have cogently argued that much of Bolzano's work on infinity is best described as the development of a theory for infinite sums. Further development of Bolzano's approach to infinity has been made by Trlifajova [17, 18]; her numerosity theory is similar to that of Benci and DiMasso [3], but differs in an important way, because it is constructive. An infinite arithmetic that is both constructive and consistent, fulfilling Hilbert's program [see 20], may be possible.

**ACKNOWLEDGEMENTS**

I would like to thank both Dr. R. de Borchgrave and Professor A. Paseau for their helpful discussion.